\documentclass[english]{svjour3}
\usepackage[T1]{fontenc}
\usepackage[latin9]{inputenc}
\usepackage{url}
\usepackage{amsmath}
\usepackage{amssymb}

\makeatletter

\newcommand{\lyxmathsym}[1]{\ifmmode\begingroup\def\b@ld{bold}
  \text{\ifx\math@version\b@ld\bfseries\fi#1}\endgroup\else#1\fi}

\numberwithin{equation}{section}
\numberwithin{figure}{section}

\makeatother

\usepackage{babel}
\begin{document}
\title{Foundations for conditional probability}
\author{Ladislav Me\v{c}í\v{r}}
\institute{Saphirion AG, Zug, Switzerland\\
lmecir@volny.cz}
\maketitle
\begin{abstract}
This article formalizes probability as a function naturally induced
by a plausible preorder on random quantities. It shows that the probability
rules including the Bayes' rule, are derivable from just this fundamental
characterization.

According to supplementary results, probability is invariably described
as a coherent function and every coherent function can be extended
to a plausibly complete function. Thus, probability can be, without
loss of generality, formalized as a plausibly complete function.

As an illustration, consider a plausibly complete probability $P$.
Then for every event $A$ and nonzero event $C$ holds that $P(A|C)=0$
if $A\wedge C=0$ and $P(A|C)=1$ if $A\wedge C=C$, no matter whether
the unconditional probability $P(C)$ is zero or whether it is defined.
In contrast to that, the common approach using the $\frac{P(A\wedge C)}{P(C)}$
ratio to define the conditional probability $P(A|C)$ leaves the probability
plausibly incomplete in general, since it leaves the conditional probability
undefined whenever $P(C)$ is zero or undefined.
\end{abstract}


\section{Introduction}

The probability foundations provided by A. N. Kolmogorov \cite{key-10}
define conditional probability as a ratio of unconditional probabilities.
A. Hájek \cite{key-5} brings several reasons why a more adequate
formalization of conditional probability is needed.

R.T. Cox \cite{key-2} contributed a theorem deriving the laws of
conditional probability from a set of postulates. According to J.
Halpern \cite{key-6}, M. J. Dupré and F. J. Tipler \cite{key-4},
J. B. Paris \cite{key-8} as well as other authors, Cox's approach
is non-rigorous. To be valid, it needs additional assumptions which
are complicated and nontrivial.

B. de Finetti \cite{key-3} developed the foundations of conditional
probability around the idea of a partially ordered algebra of random
quantities, on which existence of a real-valued function having the
fundamental properties of conditional expectation is postulated.

We take a more general approach. Instead of postulating the existence
of a real-valued function having the fundamental properties of conditional
expectation, we examine a relation on random quantities called a plausible
preorder. We show that a plausible preorder on random quantities naturally
induces a set of conditional preorders. The conditional preorders
naturally induce a conditional expectation that, in general, is a
partial function assigning elements of the extended real line to pairs
consisting of a random quantity and a nonzero event. The conditional
expectation is demonstrated to satisfy a generalized form of probability
rules. In the final section, we provide a formal description of the
notion of coherence and prove that all formalizations of probability
discussed in this article are coherent. Finally, we demonstrate that
a function is coherent if and only if it can be extended to a conditional
expectation naturally induced by a regular plausible preorder.

\section{Random quantities}

We define the notion of a random quantity using the axiomatic approach
proposed by B. de Finetti \cite{key-3}.

\subsection{Definition\label{rand-q-df}}

Let $\mathcal{T}$ denote the set of \emph{random quantities}. We
postulate that $\mathcal{T}$ is a unital associative commutative
algebra over real numbers, i.e. $\mathcal{T}$ is a set equipped with
addition, multiplication by real numbers and multiplication, such
that
\begin{itemize}
\item if $X,Y,Z\in\mathcal{T}$, then $\left(X+Y\right)+Z=X+\left(Y+Z\right)$
(associativity of addition)
\item if $X,Y\in\mathcal{T}$, then $X+Y=Y+X$ (commutativity of addition)
\item there exists an element $\mathbf{0}\in\mathcal{T}$, such that $X+\mathbf{0}=X$
for every $X\in\mathcal{T}$ (identity element of addition)
\item if $X\in\mathcal{T}$, then there exists an element $-X\in\mathcal{T}$
such that $\mathbf{0}=X+\left(-X\right)$ (inverse elements of addition)
\item if $r,s$ are real numbers and $X\in\mathcal{T}$, then $\left(rs\right)X=r\left(sX\right)$
(compatibility of multiplication by real numbers with real multiplication)
\item if $X\in\mathcal{T}$, then $1X=X$ (identity element of multiplication
by real numbers)
\item if $r$ is a real number and $X,Y\in\mathcal{T}$, then $r\left(X+Y\right)=rX+rY$
(distributivity of multiplication by real numbers with respect to
addition)
\item if $r,s$ are real numbers and $X\in\mathcal{T}$, then $\left(r+s\right)X=rX+sX$
(distributivity of multiplication by real numbers with respect to
real addition)
\item if $X,Y,Z\in\mathcal{T}$, then $\left(X.Y\right).Z=X.\left(Y.Z\right)$
(associativity of multiplication)
\item if $X,Y\in\mathcal{T}$, then $X.Y=Y.X$ (commutativity of multiplication)
\item if $X,Y,Z\in\mathcal{T}$, then $\left(X+Y\right).Z=X.Z+Y.Z$ (distributivity)
\item if $r,s$ are real numbers and $X,Y\in\mathcal{T}$, then $\left(rX\right).\left(sY\right)=\left(rs\right)\left(X.Y\right)$
(compatibility with multiplication by real numbers)
\item there exists an element $\mathbf{1}\in\mathcal{T}$, such that $X.\mathbf{1}=X$
for every $X\in\mathcal{T}$ (identity element of multiplication)
\end{itemize}

\subsection{Canonical embedding of real numbers}

Per \ref{rand-q-df}, $\mathcal{T}$ has got an identity element of
multiplication that we can denote $\mathbf{1}$. We define a function
$F$ from the set of real numbers $\mathbb{R}$ to $\mathcal{T}$
such that $F(r)=r\mathbf{1}$ for every real number $r$. Defined
this way, $F$ is a map embedding the set of real numbers in $\mathcal{T}$.
We call this embedding the \emph{canonical embedding of real numbers}
in $\mathcal{T}$. Using the canonical embedding of real numbers,
instead of writing $r\mathbf{1}\in\mathcal{T}$ for a real $r$, we
simply write $r\in\mathcal{T}$ from now on.

\subsection{Motivational example\label{rand-q-mot}}

Alice is going to throw a coin in the presence of a notary. Bob knows
that Carol shall pay him a specific amount $S_{H}$ if Alice throws
heads, and a specific amount $S_{T}$ if Alice throws tails.

Bob conceives a set $\mathcal{T}$ containing pairs of real numbers
$(X_{H},X_{T})$ and defines
\begin{itemize}
\item addition: $(X_{H},X_{T})+(Y_{H},Y_{T})=(X_{H}+Y_{H},X_{T}+Y_{T})$
\item multiplication by real numbers: $r(X_{H},X_{T})=(rX_{H},rX_{T})$
\item multiplication: $(X_{H},X_{T}).(Y_{H},Y_{T})=(X_{H}Y_{H},X_{T}Y_{T})$
\end{itemize}
Bob's $\mathcal{T}$ with these operations is unital, since $(1,1)$
is its identity element of multiplication, associative and commutative
algebra over reals. Per \ref{rand-q-df}, the elements of $\mathcal{T}$
are random quantities. Denoting $H=(1,0)$ and $T=(0,1)$ and using
the canonical embedding of real numbers, $H+T=1$ and $H.T=0$. Carol's
payment is represented by random quantity $S=(S_{H},S_{T})=S_{H}H+S_{T}T$.

\section{Events}

\subsection{Definition\label{events-df}}

Let $A$ be a random quantity. We say that $A$ is an \emph{event}
if it is idempotent, i.e. if $A.A=A$. We denote the set of events
$\mathcal{E}(\mathcal{T})$ and the set of nonzero events $\mathcal{E}_{0}(\mathcal{T})$.
On $\mathcal{E}(\mathcal{T})$ we define
\begin{itemize}
\item negation, for event $A$ its negation $\neg A$ is defined as $1-A$
\item conjunction, for events $A,B$ their conjunction $A\wedge B$ is defined
as $A.B$
\item disjunction, for events $A,B$ their disjunction $A\vee B$ is defined
as $A+B-A.B$
\end{itemize}
With these operations,
\begin{itemize}
\item $1$ is the identity element of conjunction, i.e. if $A$ is an event,
then $A\wedge1=1\wedge A=A$
\item $0$ is the identity element of disjunction, i.e. if $A$ is an event,
then $A\vee0=0\vee A=A$ and
\item $\mathcal{E}(\mathcal{T})$ is a Boolean algebra.
\end{itemize}

\subsection{Natural order}

Let $\mathcal{A}$ be a Boolean algebra and $A,B$ be its elements.
The natural order on $\mathcal{A}$ is defined so that $A\leq B$
if $A\wedge B=A$.

Then
\begin{itemize}
\item the natural order is a partial order,
\item the minimal element of $\mathcal{A}$ in the natural order is 0 and
\item the maximal element of $\mathcal{A}$ in the natural order is 1.
\end{itemize}
Since $\mathcal{E}(\mathcal{T})$ is a Boolean algebra, there is a
natural order on $\mathcal{E}(\mathcal{T})$.

\subsection{Atoms}

Let $\mathcal{A}$ be a Boolean algebra. We say that $D$ is an \emph{atom}
of $\mathcal{A}$, if $D\in\mathcal{A}$ and $D$ is a minimal nonzero
element of $\mathcal{A}$ in the natural order.

We say that $\mathcal{A}$ is \emph{atomic} if for every nonzero element
$A\in\mathcal{A}$ there is an atom $D\in\mathcal{A}$ such that $D\le A$
in the natural order.

\subsection{Motivational example}

Consider the algebra $\mathcal{T}$ defined in \ref{rand-q-mot}.
Then
\begin{itemize}
\item $\mathcal{E}\left(\mathcal{T}\right)=\{0,H,T,1\}$
\item $\mathcal{E}_{0}\left(\mathcal{T}\right)=\{H,T,1\}$
\item the atoms of $\mathcal{E}\left(\mathcal{T}\right)$ are $H$ and $T$
and
\item $\mathcal{E}\left(\mathcal{T}\right)$ is atomic.
\end{itemize}

\subsection{Positive combinations of nonzero events\label{pos-comb-nz-events}}

Let $n\ge1$, $p_{1},\ldots,p_{n}$ be positive real numbers and $C_{1},\ldots,C_{n}$
be nonzero events. Then $0\ne\sum_{i=1}^{n}p_{i}C_{i}$.

\section{Plausible preorder}

\subsection{Definition\label{pl-preord-df}}

We say that a relation $\lesssim$ on $\mathcal{T}$ is a \emph{plausible
preorder} if it has these properties:

\subsubsection{Plausible property\label{plausible-preord-plausible-property}}

If $A$ is an event, then $0\lesssim A$.

\subsubsection{Additive property\label{pl-preord-Additive-property}}

If $0\lesssim X$ and $0\lesssim Y$, then $0\lesssim X+Y$.

\subsubsection{Multiplicative property\label{plaus-preor-multiplicative-property}}

If $0\lesssim X$ and $q$ is a nonnegative real number, then $0\lesssim qX$.

\subsubsection{Extension property\label{pl-preord-Extension-property}}

$X\lesssim Y$ if and only if $0\lesssim Y-X$.

\subsection{Motivational example\label{pl-preord-mot}}

In the algebra $\mathcal{T}$ described in \ref{rand-q-mot}, Bob
defines
\begin{itemize}
\item $0\lesssim X=(X_{H},X_{T})$ if $0\leq X_{H}+X_{T}$
\item $X\lesssim Y$ if $0\lesssim Y-X$
\end{itemize}
It is easy to verify that $\lesssim$ is a plausible preorder.

\subsection{The greatest plausible preorder}

Relation $\lesssim=\mathcal{T}\times\mathcal{T}$ is a plausible preorder
and the greatest relation\emph{ }on $\mathcal{T}$ with respect to
inclusion.

\subsection{Properties\label{pl-preord-properties}}

\subsubsection{Reflexivity}

A plausible preorder is reflexive.

\subsubsection{Transitivity}

A plausible preorder is transitive.

\subsubsection{Relation to the natural order of events.\label{preorder-to-natural-order}}

A plausible preorder contains the natural order of events as its subset.

\subsubsection{Relation to the order of real numbers.\label{preorder-to-order-of-reals}}

A plausible preorder contains the order of real numbers as its subset.

\subsubsection{Intersection of a set of plausible preorders}

If $\mathcal{X}$ is a nonempty set containing plausible preorders,
then $\bigcap\mathcal{X}$ is a plausible preorder.

\subsubsection{The smallest plausible preorder containing a relation}

If $R$ is a relation on $\mathcal{T}$, then there is a relation
$\lesssim$ that is the smallest plausible preorder with respect to
inclusion containing $R$.

\subsubsection{Subadditivity\label{pl-preord-subadditivity}}

If $A_{1},\ldots,A_{n}$ are events and $\lesssim$ is a plausible
preorder, then
\[
\bigvee_{i=1}^{n}A_{i}\lesssim\sum_{i=1}^{n}A_{i}
\]

\section{Plausible equivalence}

\subsection{Definition}

Let $\lesssim$ be a plausible preorder. We say that a relation $\sim$
is the \emph{equivalence part} of $\lesssim$ if for any $X,Y$ holds
that $X\sim Y$ if and only if $(X\lesssim Y)\wedge(Y\lesssim X)$.
We say that a relation $\sim$ on $\mathcal{T}$ is a \emph{plausible
equivalence} if there is a plausible preorder $\lesssim$ such that
$\sim$ is its equivalence part.

\subsection{Motivational example}

Let $\lesssim$ be the plausible preorder defined in \ref{pl-preord-mot}.
Then $0\sim X=\left(X_{H},X_{T}\right)$ if and only if $0=X_{H}+X_{T}$. 

\subsection{Fundamental properties}

\subsubsection{Plausible property\label{pl-equiv-pl-property}}

Let $p_{1},\ldots,p_{n}$ be positive real numbers, $A_{1},\ldots,A_{n}$
be events and $0\sim\sum_{i=1}^{n}p_{i}A_{i}$. Then $0\sim A_{i}$
for every $i\in\left\{ 1,\ldots,n\right\} $.

\subsubsection{Reflexivity}

$0\sim0$.

\subsubsection{Additive property}

If $0\sim X$ and $0\sim Y$, then $0\sim X+Y$.

\subsubsection{Multiplicative property}

If $0\sim X$ and $r$ is a real number, then $0\sim rX$.

\subsubsection{Extension property}

$X\sim Y$ if and only if $0\sim Y-X$.

\subsection{Sufficiency of the fundamental properties\label{pl-equiv-sufficiency}}

Every relation $\sim$ having the fundamental properties of a plausible
equivalence is a plausible equivalence.

\section{Plausible strict partial order}

\subsection{Definition\label{strict-p-o}}

Let $\lesssim$ be a plausible preorder. We say that a relation $\lnsim$
is the \emph{strict part} of $\lesssim$ if for any $X,Y$ holds that
$X\lnsim Y$ if and only if $\left(X\lesssim Y\right)\land\lnot\left(Y\lesssim X\right)$.
We also say that a relation $\lnsim$ is a \emph{plausible strict
partial order} if there is a plausible preorder $\lesssim$ such that
$\lnsim$ is its strict part.

\subsection{Motivational example}

Let $\lesssim$ be the plausible preorder defined in \ref{pl-preord-mot}.
Then $0\lnsim X=\left(X_{H},X_{T}\right)$ if and only if $0<X_{H}+X_{T}$. 

\subsection{Fundamental properties\label{strict-p-o-fund}}

\subsubsection{Plausible property\label{strict-p-o-pl-property}}

If $A$ is an event, $\neg(0\lnsim A)$ and $0\lnsim X$, then $0\lnsim X+A$
and $0\lnsim X-A$.

\subsubsection{Antireflexivity\label{strict-p-o-Antireflexivity}}

$\neg(0\lnsim0)$.

\subsubsection{Additive property}

If $0\lnsim X$ and $0\lnsim Y$, then $0\lnsim X+Y$.

\subsubsection{Multiplicative property}

If $0\lnsim X$ and $p$ is a positive real number, then $0\lnsim pX$.

\subsubsection{Extension property}

$X\lnsim Y$ if and only if $0\lnsim Y-X$.

\subsection{Sufficiency of the fundamental properties\label{strict-p-o-sufficiency}}

Every relation $\lnsim$ having the fundamental properties of a plausible
strict partial order is a plausible strict partial order.

\section{Conditional preorder}

\subsection{Definition\label{Conditional-preorder-df}}

Let $\lesssim$ be a plausible preorder and $C$ be an event. We define
the \emph{conditional preorder} $\lesssim_{C}$ so that $X\lesssim_{C}Y$
if $X.C\lesssim Y.C$.

\subsection{Properties}
\begin{itemize}
\item a conditional preorder is a plausible preorder
\item $\lesssim_{1}$ is identical with $\lesssim$
\item $\lesssim_{0}$ is the greatest plausible preorder
\item $\lnsim_{0}$ is empty, i.e. there are no random quantities $X,Y$
such that $X\lnsim_{0}Y$
\end{itemize}

\section{Regularity of a plausible preorder}

\subsection{Definition}

We say that a plausible preorder $\lesssim$ is
\begin{itemize}
\item \emph{degenerate} if $0\sim1$
\item \emph{regular} if for every nonzero event $C$ holds $0\lnsim C$
\end{itemize}

\subsection{Motivational example}

The plausible preorder $\lesssim$ defined in \ref{pl-preord-mot}
is regular.

\subsection{Properties}
\begin{itemize}
\item A conditional preorder $\lesssim_{C}$ is degenerate if and only if
$0\sim C$.
\item The greatest plausible preorder is degenerate.
\item If a plausible preorder is degenerate, then for every pair of random
quantities $X,Y$ in the linear span of $\mathcal{E}(\mathcal{T})$
holds $X\sim Y$. In particular, for every pair of real numbers $r,s$
holds $r\sim s$, and for every pair of events $A,B$ holds $A\sim B$.
\item A plausible preorder is nondegenerate if and only if it coincides
with the order of real numbers on $\mathbb{R}$.
\end{itemize}

\section{Extended real line}

\subsection{Definition}

We define the \emph{extended real line} as the set $\overline{\mathbb{R}}=\mathbb{R}\cup\{-\infty,+\infty\}$,
where $\mathbb{R}$ is the set of real numbers.

\subsection{Order}

We extend the order of real numbers to $\overline{\mathbb{R}}$ so
that $-\infty\le x\le+\infty$ for every $x\in\overline{\mathbb{R}}$,
turning $\overline{\mathbb{R}}$ into a linearly ordered set. In this
order, every subset $U$ of $\overline{\mathbb{R}}$ has both the
least upper bound (\emph{supremum}) denoted $\sup U$ and the greatest
lower bound (\emph{infimum}) denoted $\inf U$. In particular,
\begin{itemize}
\item $\sup\emptyset=-\infty$ 
\item $\inf\emptyset=+\infty$ 
\item $\sup\mathbb{R}=\sup\overline{\mathbb{R}}=+\infty$ 
\item $\inf\mathbb{R}=\inf\overline{\mathbb{R}}=-\infty$ 
\end{itemize}

\subsection{Arithmetic}

We extend the artithmetic operations on real numbers to $\overline{\mathbb{R}}$
so that
\begin{itemize}
\item if $x\ne-\infty$, then $(+\infty)+x=x+(+\infty)=+\infty$
\item if $x\ne+\infty$, then $(-\infty)+x=x+(-\infty)=-\infty$
\item if $x>0$, then $x.(+\infty)=(+\infty).x=+\infty$
\item if $x>0$, then $x.(-\infty)=(-\infty).x=-\infty$
\item if $x<0$, then $x.(+\infty)=(+\infty).x=-\infty$
\item if $x<0$, then $x.(-\infty)=(-\infty).x=+\infty$
\item if $x$ is a real number, then $\frac{x}{+\infty}=\frac{x}{-\infty}=0$
\item if $x$ is a positive real number, then $\frac{+\infty}{x}=+\infty$
\item if $x$ is a positive real number, then $\frac{-\infty}{x}=-\infty$
\item if $x$ is a negative real number, then $\frac{+\infty}{x}=-\infty$
\item if $x$ is a negative real number, then $\frac{-\infty}{x}=+\infty$
\end{itemize}
Other expressions than the above are undefined. For example, the expressions
\begin{itemize}
\item $(+\infty)+(-\infty)$
\item $(-\infty)+(+\infty)$
\item $0.(+\infty)$
\item $0.(-\infty)$
\item $(+\infty).0$
\item $(-\infty).0$
\item $\frac{x}{0}$
\item $\frac{+\infty}{+\infty}$
\item $\frac{+\infty}{-\infty}$
\item $\frac{-\infty}{+\infty}$
\item $\frac{-\infty}{-\infty}$
\end{itemize}
are all undefined.

\section{Expectation naturally induced by a plausible preorder}

\subsection{Definition\label{Expectation}}

Let $\lesssim$ be a plausible preorder and $X$ be a random quantity.
Denoted $E\left(X\right)$, the \emph{expectation of} $X$ (more precisely,
the \emph{expectation of} $X$ \emph{naturally induced by} $\lesssim$)
is
\begin{itemize}
\item a real number $x$, if for every positive real number $\epsilon$
holds $-\epsilon\lnsim X-x\lnsim\epsilon$
\item $+\infty$, if for every real number $y$ holds $y\lnsim X$
\item $-\infty$, if for every real number $y$ holds $X\lnsim y$
\item not defined, if none of the above holds
\end{itemize}

\subsection{Motivational example\label{expectation-mot}}

Consider the plausible preorder defined in example \ref{pl-preord-mot}.
For every random quantity $X=(X_{H},X_{T})$ holds that $E(X)=\frac{1}{2}X_{H}+\frac{1}{2}X_{T}$.
In particular, $E(H)=E(T)=\frac{1}{2}$.

\subsection{Relation to regularity of plausible preorder\label{expect-to-regul}}

Let $\lesssim$ be a plausible preorder and $r$ be a real number.
\begin{itemize}
\item if $0\lnsim1$, then $E(r)=r$
\item if $0\sim1$, then $E(r)$ is not defined
\item if $\lesssim$ is the maximal plausible preorder and $X$ is a random
quantity, then $E\left(X\right)$ is not defined
\end{itemize}

\subsection{Preorder consistency\label{Preorder-consistency}}

Let $r$ be a real number and $X,Y$ be random quantities. \ref{Expectation}
implies that
\begin{itemize}
\item if $E(X)$ exists and $r<E(X)$, then $r\lnsim X$
\item if $E(X)$ exists and $E(X)<r$, then $X\lnsim r$ 
\item if both $E(X)$ and $E(Y)$ exist and $E(X)<E(Y)$, then $X\lnsim Y$
\item if both $E(X)$ and $E(Y)$ exist and $X\lesssim Y$, then $E(X)\leq E(Y)$
\end{itemize}

\subsection{Existence and uniqueness\label{Expectation-existence}}

Let $\lesssim$ be a plausible preorder and $X$ be a random quantity.
Then $X$ has expectation if and only if in the extended real line
$\overline{\mathbb{R}}$ holds that $\sup\{r\in\mathbb{R}|r\lnsim X\}=\inf\{r\in\mathbb{R}|X\lnsim r\}$.
In such case, $E(X)=\sup\{r\in\mathbb{R}|r\lnsim X\}=\inf\{r\in\mathbb{R}|X\lnsim r\}$.

\section{Conditional expectation\label{Conditional-expectation}}

\subsection{Definition\label{cond-expect-df}}

Let $\lesssim$ be a plausible preorder, $X$ be a random quantity
and $C$ be an event. Denoted $E(X|C)$, the \emph{conditional expectation
of} $X$ \emph{given} $C$ is defined using \ref{Expectation} as
the expectation of $X$ naturally induced by the conditional preorder
$\lesssim_{C}$.

\subsection{Motivational example\label{cond-expect-mot}}

Consider the plausible preorder defined in \ref{pl-preord-mot}. For
every random quantity $X=(X_{H},X_{T})$ holds that $E(X|H)=X_{H}$
and $E(X|T)=X_{T}$.

\subsubsection{Relation to regularity of plausible preorder\label{conditional-reg}}

Let $A,C$ be events and $r$ be a real number. Then
\begin{itemize}
\item if $0\lnsim C$ and $r$ is a real number, then $E\left(r|C\right)=E(rC|C)=r$
\item if $0\sim C$ then neither $E(r|C)$ nor $E(rC|C)$ is defined
\end{itemize}

\subsection{As a function}

Per \ref{expect-to-regul} and \ref{Expectation-existence}, conditional
expectation is a partial function from $\mathcal{T}\times\mathcal{E}_{0}(\mathcal{T})$
to $\overline{\mathbb{R}}$.

\section{Rules\label{expectation-rules}}

Let $\lesssim$ be a plausible preorder, $E(X|C)$ be the conditional
expectation naturally induced by $\lesssim$, $X,Y$ be random quantities,
$B,C,D$ be events and $r$ be a real number. The rules the conditional
expectation follows are:

\subsection{Consistency\label{Consistency}}

$E(X|C)$ exists if and only if $E(X.C|C)$ exists. In case it exists,
\[
E(X|C)=E(X.C|C)
\]

\subsection{Real additivity}

If $E(X|C)$ exists, then
\[
E(X+r|C)=E(X|C)+r
\]

\subsection{General additivity}

If the expression $E(X|C)+E(Y|C)$ makes sense, then
\[
E(X+Y|C)=E(X|C)+E(Y|C)
\]

\subsection{Homogeneity\label{Homogeneity}}

If the expression $rE(X|C)$ makes sense, then
\[
E(rX|C)=rE(X|C)
\]

\subsection{Conditional probability\label{Conditional-probability}}

If $E(C|D)$ exists, we, compatibly with Thomas Bayes \cite{key-1},
denote
\[
P(C|D)=E(C|D)
\]
and say that it is the \emph{conditional probability of $C$ given
$D$}. For $P(C|1)$ we also use a simpler notation $P(C)$.

\subsection{Monotonicity}

If both $P(B|D)$ and $P(C|D)$ exist and $B\leq C$ in the natural
order of events, then
\[
P(B|D)\leq P(C|D)
\]

\subsection{Minimal and maximal probability\label{Minimal-and-maximal-pr}}

If $P(C|D)$ exists, then
\[
0=P(0|D)\leq P(C|D)\leq P(D|D)=P(1|D)=1
\]

\subsection{Completeness}
\begin{itemize}
\item If $P(B|D)=0$ and $C\leq B$ in the natural order of events, then
\[
P(C|D)=0
\]
\item If $P(B|D)=1$ and $B\leq C$ in the natural order of events, then
\[
P(C|D)=1
\]
\end{itemize}

\subsection{Subadditivity}

If $A_{1},\ldots,A_{n}$ are events and all of $P(A_{1}|D),\ldots,P(A_{n}|D),P(\bigvee_{i=1}^{n}A_{i}|D)$
exist, then 
\[
P(\bigvee_{i=1}^{n}A_{i}|D)\leq\sum_{i=1}^{n}P(A_{i}|D)
\]

\subsection{Bayes' rule\label{subsec:Bayes'-rule}}

\subsubsection{Chain form\label{Chain}}

If the expression $E(X|C.D).P(C|D)$ makes sense, then
\[
E(X.C|D)=E(X|C.D).P(C|D)
\]

\subsubsection{Chain form, zero \emph{E(X|C.D)\label{Chain-zero-e}}}

If $E(X|C.D)=0$, then
\[
E(X.C|D)=0
\]

\subsubsection{Chain form, zero \emph{P(C|D)\label{Chain-zero-p}}}

If $P\left(C|D\right)=0$ and there is a real number $p$ such that
$-p\lesssim_{C.D}X\lesssim_{C.D}p$, then
\[
E\left(X.C|D\right)=0
\]

\subsubsection{Chain form, infinite \emph{E(X|C.D)\label{Chain-infinite-e}}}

If $E\left(X|C.D\right)\in\{-\infty,+\infty\}$ and there is a positive
real number $p$ such that $p\lesssim_{D}C$, then
\[
E\left(X.C|D\right)=E\left(X|C.D\right)
\]

\subsubsection{Conditional form\label{Cond-form}}

If the expression $\frac{E(X.C|D)}{P(C|D)}$ makes sense, then
\[
E(X|C.D)=\frac{E(X.C|D)}{P(C|D)}
\]

\subsubsection{Conditional form, zero \emph{E(X.C|D)\label{Cond-form-zero-e}}}

If $E\left(X.C|D\right)=0$ and there is a positive real number $p$
such that $p\lesssim_{D}C$, then 
\[
E\left(X|C.D\right)=0
\]

\subsubsection{Conditional form, infinite \emph{E(X.C|D)\label{Cond-form-inf-e}}}

If $E(X.C|D)\in\{-\infty,+\infty\}$, then
\[
E(X|C.D)=E(X.C|D)
\]

\subsubsection{Conditional form, zero \emph{P(C|D)\label{Cond-form-zero-p}}}

If $P(C|D)=0$ and there is a real number $p$ such that $1\lesssim_{D}pX.C$,
then $p\neq0$ and 
\[
E\left(X|C.D\right)=\frac{+\infty}{p}
\]

\subsubsection{P form\label{P-form}}

If the expression $\frac{E(X.C|D)}{E(X|C.D)}$ makes sense, then
\[
P(C|D)=\frac{E(X.C|D)}{E(X|C.D)}
\]

\subsubsection{P form, zero \emph{E(X.C|D)\label{P-form-zero-e(xc|d)}}}

If $E\left(X.C|D\right)=0$ and there is a real number $p$ such that
$1\lesssim_{C.D}pX$, then 
\[
P\left(C|D\right)=0
\]

\subsubsection{P form, infinite \emph{E(X|C.D)\label{P-form-inf-e(x|cd)}}}

If $E\left(X|C.D\right)\in\left\{ -\infty,+\infty\right\} $ and there
is a real number $p$ such that $-p\lesssim_{D}X.C\lesssim_{D}p$,
then 
\[
P\left(C|D\right)=0
\]

\section{Coherence\label{coherence}}

\subsection{Definition\label{Coherence-df}}

We say that $PV$ \emph{is a coherent function} if it is a partial
function from $\mathcal{T}\times\mathcal{E}\left(\mathcal{T}\right)$
to $\overline{\mathbb{R}}$ such that if
\begin{itemize}
\item $n\geq0,m\geq1$ are integers
\item $q_{1},\ldots,q_{n}$ are nonnegative real numbers
\item $r_{1},\ldots,r_{m}$, $s_{1},\ldots,s_{m}$ are real numbers
\item $C_{1},\ldots,C_{n}$, $D_{1},\ldots,D_{m}$ are events
\item $X_{1},\ldots,X_{m}$ are random quantities
\item $r_{j}(PV(X_{j}|D_{j})+s_{j})>0$ for every $j\in\left\{ 1,\ldots,m\right\} $
\end{itemize}
then 
\[
0\neq\sum_{i=1}^{n}q_{i}C_{i}+\sum_{j=1}^{m}r_{j}(X_{j}+s_{j}).D_{j}
\]

\subsection{Kolmogorovian plausible values\label{Kolmogorovian-pv}}

Let $\mathbb{R}_{\geq0}$ denote the set of nonnegative real numbers.
We say that $PV$ is a \emph{Kolmogorovian plausible value} if
\begin{itemize}
\item $PV$ is a function from $\mathcal{F}$ to $\mathbb{R}_{\geq0}$,
where $\mathcal{F}$ is a nonempty subset of $\mathcal{E}(\mathcal{T})$
closed under negation and conjunction
\item $PV(1)=1$ (unitarity)
\item if $A,B\in\mathcal{F}$ and $A.B=0$, then $PV(A+B)=PV(A)+PV(B)$
(additivity)
\end{itemize}
Using the notation $PV(A|1)=PV(A)$, we can handle every Kolmogorovian
plausible value as a function from $\mathcal{F}\times\{1\}$ to $\mathbb{R}_{\geq0}$.

\subsubsection{Coherence\label{Kolmogorovian-coherence}}

Every Kolmogorovian plausible value is coherent.

\subsection{Coxian plausible values\label{Coxian-pv}}

We say that $PV$ is a \emph{Coxian plausible value} if
\begin{itemize}
\item $PV$ is a function from $\mathcal{F}\times\mathcal{F}_{0}$ to $\mathbb{R}_{\geq0}$,
where $\mathcal{F}$ is a nonempty subset of $\mathcal{E}(\mathcal{T})$
closed under negation and conjunction and $\mathcal{F}_{0}$ is the
set containing all elements of $\mathcal{F}$ except for $0$
\item if $C\in\mathcal{F}_{0}$, then $PV(C|C)>0$ (positivity)
\item if $A\in\mathcal{F}$ and $C\in\mathcal{F}_{0}$, then $PV(1-A|C)=1-PV(A|C)$
(negation formula)
\item if $A,C,D$ are elements of $\mathcal{F}$ and $C.D\neq0$, then $PV(A.C|D)=PV(A|C.D).PV(C|D)$
(Bayes' rule)
\end{itemize}

\subsubsection{Basic properties}

Let $PV$ be a function from $\mathcal{F}\times\mathcal{F}_{0}$ to
$\mathbb{R}_{\geq0}$ that satisfies definition \ref{Coxian-pv},
$A\in\mathcal{F}$ and $C\in\mathcal{F}_{0}$. Then
\begin{itemize}
\item $PV(C|C)=1$
\item $PV(1|C)=1$
\item $PV(0|C)=0$
\item $PV(A.C|C)=PV(A|C)$
\end{itemize}

\subsubsection{Sum rule\label{Coxian-sum-rule}}

Let $PV$ be a function from $\mathcal{F}\times\mathcal{F}_{0}$ to
$\mathbb{R}_{\geq0}$ that satisfies definition \ref{Coxian-pv}.
Let $A,B\in\mathcal{F}$ such that $A.B=0$ and let $C\in\mathcal{F}_{0}$.
Then 
\[
PV(A+B|C)=PV(A|C)+PV(B|C)
\]

\subsubsection{Subadditivity\label{Coxian-subadditivity}}

Let $PV$ be a function from $\mathcal{F}\times\mathcal{F}_{0}$ to
$\mathbb{R}_{\geq0}$ that satisfies definition \ref{Coxian-pv}.
Let $A_{1},\ldots,A_{n}\in\mathcal{F}$ and let $C\in\mathcal{F}_{0}$.
Then 
\[
PV(\bigvee_{i=1}^{n}A_{i}|C)\leq\sum_{i=1}^{n}PV(A_{i}|C)
\]

\subsubsection{Coherence\label{Coxian-coherence}}

Every Coxian plausible value is coherent.

\subsection{Dupré-Tiplerian plausible values\label{Dupre-Tiplerian-pv}}

We say that $PV$ is a \emph{Dupré-Tiplerian plausible value} if
\begin{itemize}
\item $PV$ is a partial function from $\mathcal{T}\times\mathcal{C}$ to
$\mathbb{R}$, where $\mathcal{C}$ is a subset of $\mathcal{E}_{0}(\mathcal{T})$
closed under disjunction
\item if $A$ is an event and $C\in\mathcal{C}$, then $PV(A|C)$ exists
and $PV(A|C)\geq0$ (nonnegativity)
\item if $C\in\mathcal{C}$, then $PV(C|C)>0$ (positivity)
\item if $r\in\mathbb{R}$, $C\in\mathcal{C}$, $X$ is a random quantity
and $PV(X|C)$ exists, then $PV(rX|C)=r.PV(X|C)$ (homogeneity)
\item if $C\in\mathcal{C}$, $X$, $Y$ are random quantities and both $PV(X|C)$
and $PV(Y|C)$ exist, then $PV(X+Y|C)=PV(X|C)+PV(Y|C)$ (additivity)
\item if $C.D\in\mathcal{C}$, $D\in\mathcal{C}$, $X$ is a random quantity
and $PV(X|C.D)$ exists, then $PV(X.C|D)=PV(X|C.D).PV(C|D)$ (Bayes'
rule)
\end{itemize}

\subsubsection{Subadditivity\label{Dupre-Tiplerian-subadditivity}}

Let $PV$ be a function from $\mathcal{T}\times\mathcal{C}$ to $\mathbb{R}$
that satisfies definition \ref{Dupre-Tiplerian-pv}. Let $A_{1},\ldots,A_{n}$
be events and let $C\in\mathcal{C}$. Then 
\[
PV(\bigvee_{i=1}^{n}A_{i}|C)\leq\sum_{i=1}^{n}PV(A_{i}|C)
\]

\subsubsection{Coherence\label{Dupre-Tiplerian-coherence}}

Every Dupré-Tiplerian plausible value is coherent.

\subsection{Characterizations\label{Coherent functions}}

Let $PV$ be a partial function from $\mathcal{T}\times\mathcal{E}(\mathcal{T})$
to $\overline{\mathbb{R}}$. Then the following characterizations
are equivalent:
\begin{enumerate}
\item $PV$ is coherent
\item $PV$ can be extended to conditional expectation naturally induced
by a regular plausible preorder
\item $PV$ can be extended to conditional expectation naturally induced
by a plausible preorder
\end{enumerate}

\subsection{Probability as a plausibly complete function}

We say that a function is \emph{plausibly complete}, if it is naturally
induced by a regular plausible preorder. According to \ref{Kolmogorovian-coherence},
\ref{Coxian-coherence}, \ref{Dupre-Tiplerian-coherence} and \ref{Coherent functions},
we can, without loss of generality, characterize probability as a
plausibly complete function.

\section{Conclusion}

Assigning the role of a primitive notion to the notion of a plausible
preorder, our formalization offers a different perspective on the
foundations of probability than the formalizations discussed in the
introduction. Our approach neither forces us to define conditional
probability by a ratio of unconditional probabilities which is criticized
as inadequate, nor does it force us to postulate conditional probability
to have other properties open to doubt. The formalization is supported
by theorem \ref{Coherent functions}, confirming that it encompasses
all coherent instances of probability. We supplement it by verifying
that according to all formalizations of the probability notion discussed
in the introduction, probability is coherent. To illustrate that our
formalization satisfies the main Hájek's \cite{key-5} requirements,
consider a nonzero event $C$ such that $P\left(C\right)$ is either
zero or undefined. Because of that, the ratio $\frac{P\left(A\wedge C\right)}{P\left(C\right)}$
leaves the conditional probability $P\left(A|C\right)$ undefined.
On the other hand, once probability is coherent, theorem \ref{Coherent functions}
confirms that it can be extended to a plausibly complete function,
i.e. to a conditional expectation naturally induced by a regular plausible
preorder. Definition \ref{cond-expect-df} applied to a regular plausible
preorder yields that $P\left(A|C\right)=0$ if $A\wedge C=0$ and
$P\left(A|C\right)=1$ if $A\wedge C=C$, no matter whether $P\left(C\right)$
is zero or whether it is defined.

\section{Appendix}

\subsection{Proof of \ref{pos-comb-nz-events}}

Let $n\geq1$, $p_{1},\ldots,p_{n}$ be positive real numbers and
$C_{1},\ldots,C_{n}$ be nonzero events. Let $\mathcal{A}$ be the
Boolean subalgebra of $\mathcal{E}\left(\mathcal{T}\right)$ generated
by $C_{1},\ldots,C_{n}$. Since $\mathcal{A}$ is finitely generated,
it is finite and atomic \cite{key-9}. Since $C_{1}$ is nonzero,
there is an atom $D$ of $\mathcal{A}$ such that $D\le C_{1}$ in
the natural order, i.e. $C_{1}.D=D$. Since $D$ is an atom, for every
$i\in\{2,\ldots,n\}$ either $C_{i}.D=D$ or $C_{i}.D=0$. Therefore,
$\left(\sum_{i=1}^{n}p_{i}C_{i}\right).D=p_{1}D+\sum_{i=2}^{n}p_{i}C_{i}.D=pD$,
where $p\geq p_{1}$ is a real number. Since $D$ is nonzero and $p$
is a positive real number, $pD$ is nonzero. Since $\left(\sum_{i=1}^{n}p_{i}C_{i}\right).D$
is nonzero, $\sum_{i=1}^{n}p_{i}C_{i}$ is nonzero.

\subsection{Proof of \ref{pl-preord-subadditivity}}

Let $A_{1},\ldots,A_{n}$ be events and $\lesssim$ be a plausible
preorder. Due to reflexivity of $\lesssim$, the inequality $\bigvee_{i=1}^{n}A_{i}\lesssim\sum_{i=1}^{n}A_{i}$
holds for $n=0$ and $n=1$. Let the inequality $\bigvee_{i=1}^{n}A_{i}\lesssim\sum_{i=1}^{n}A_{i}$
hold for some integer $n\geq1$ and arbitrary events $A_{1},\ldots,A_{n}$.
Let $A_{1},\ldots,A_{n+1}$ be events. Define events $B_{1}=A_{1},\ldots,B_{n-1}=A_{n-1},B_{n}=A_{n}\vee A_{n+1}$.
Then $\bigvee_{i=1}^{n+1}A_{i}=\bigvee_{i=1}^{n}B_{i}\lesssim\sum_{i=1}^{n}B_{i}=\sum_{i=1}^{n+1}A_{i}-A_{n}.A_{n+1}\lesssim\sum_{i=1}^{n+1}A_{i}$.
By mathematical induction the inequality holds for every integer $n$.

\subsection{Proof of \ref{pl-equiv-pl-property}}

Let $\sim$ be the equivalence part of a plausible equivalence $\lesssim$,
let $p_{1},\ldots,p_{n}$ be positive real numbers, $A_{1},\ldots,A_{n}$
be events and $0\sim\sum_{i=1}^{n}p_{i}A_{i}$. Without loss of generality,
we assume that $n\ge1$ and prove that $0\sim A_{1}$. Per the assumption,
$\sum_{i=1}^{n}p_{i}A_{i}\lesssim0$. Per \ref{pl-preord-df}, $-\sum_{i=2}^{n}p_{i}A_{i}\lesssim0$.
Therefore, $A_{1}=\frac{1}{p_{1}}(\sum_{i=1}^{n}p_{i}A_{i}-\sum_{i=2}^{n}p_{i}A_{i})\lesssim0$.
Together with the plausible property of $\lesssim$ guaranteeing that
$0\lesssim A_{1}$, we obtain that $0\sim A_{1}$.

\subsection{Proof of \ref{pl-equiv-sufficiency}}

Let $\sim$ be a relation having the fundamental properties of a plausible
equivalence. We define a relation $\lesssim$ so that $0\lesssim X$
if $X=U+\sum_{i=1}^{n}q_{i}A_{i}$ for some random quantity $U\sim0$,
integer $n\geq0$, nonnegative real numbers $q_{1},\ldots,q_{n}$
and events $A_{1},\ldots,A_{n}$. We also define that $X\lesssim Y$
if $0\lesssim Y-X$.

The task to verify that the relation $\lesssim$ defined this way
is a plausible preorder and that $\sim$ is its equivalence part,
is left as an exercise to the reader. Note also that the relation
$\lesssim$ defined this way is the smallest plausible preorder with
respect to inclusion such that $\sim$ is its equivalence part.

\subsection{Proof of \ref{strict-p-o-sufficiency}}

Let $\lnsim$ be a relation having the fundamental properties of a
plausible strict partial order. Let $X$ be a random quantity. We
define that $0\lesssim X$ if either $0\lnsim X$ or if there are
real numbers $r_{1},r_{2},\ldots,r_{n}$ and events $A_{1},A_{2},\ldots,A_{n}$
such that $X=\sum_{i=1}^{n}r_{i}A_{i}$ and for every $i\in{1,2,\ldots,n}$
holds $\neg(0\lnsim A_{i})$. We also define that $X\lesssim Y$ if
$0\lesssim Y-X$.

The task to verify that the relation $\lesssim$ defined this way
is a plausible preorder and that $\lnsim$ is its strict part, is
left as an exercise to the reader. Note also that the relation $\lesssim$
defined this way is the smallest plausible preorder with respect to
inclusion such that $\lnsim$ is its strict part.

\subsection{Proof of \ref{expectation-mot}}

Let $X=\left(X_{H},X_{T}\right)$. Then $X-\left(\frac{1}{2}X_{H}+\frac{1}{2}X_{T}\right)=\left(\frac{1}{2}X_{H}-\frac{1}{2}X_{T}\right)H+\left(\frac{1}{2}X_{T}-\frac{1}{2}X_{H}\right)T$.
Since $\left(\frac{1}{2}X_{H}-\frac{1}{2}X_{T}\right)+\left(\frac{1}{2}X_{T}-\frac{1}{2}X_{H}\right)=0$,
\ref{pl-preord-mot} implies that $X-\left(\frac{1}{2}X_{H}+\frac{1}{2}X_{T}\right)\sim0$.
Let $\epsilon$ be a positive real number. Per \ref{pl-preord-mot},
$0\lnsim\epsilon$. Therefore, $-\epsilon\lnsim0\lesssim X-\left(\frac{1}{2}X_{H}+\frac{1}{2}X_{T}\right)\lesssim0\lnsim\epsilon$.

\subsection{Proof of \ref{Expectation-existence}}

Let $\lesssim$ be a plausible preorder and $X$ be a random quantity.

If $E(X)=+\infty$, then in the extended real line, $\sup\{r\in\mathbb{R}|r\lnsim X\}=\sup\mathbb{R}=+\infty=\inf\emptyset=\inf\{r\in\mathbb{R}|X\lnsim r\}$.

If $E(X)=-\infty$, then $\sup\{r\in\mathbb{R}|r\lnsim X\}=\sup\emptyset=-\infty=\inf\mathbb{R}=\inf\{r\in\mathbb{R}|X\lnsim r\}$.

If $E(X)=x\in\mathbb{R}$, then per \ref{Preorder-consistency}, for
$r<x$ holds that $r\lnsim X$ and for $x<r$ holds that $X\lnsim r$.
Due to antireflexivity and transitivity of $\lnsim$, there is no
real number $r$ such that both $r\lnsim X$ and $X\lnsim r$. Therefore,
if $r\lnsim X$ then also $r\le x$. Similarly, if $X\lnsim r$, then
also $x\le r$. This means that $x$ is an upper bound of $\{r\in\mathbb{R}|r\lnsim X\}$.
$x$ is also the least upper bound of $\{r\in\mathbb{R}|r\lnsim X\}$,
since for every real number $r$ smaller than $x$ there is a greater
real number $s$ that is still smaller than $x$, which means that
$r<x$ is not an upper bound of $\{r\in\mathbb{R}|r\lnsim X\}$. Similarly,
$x$ is the greatest lower bound of $\{r\in\mathbb{R}|X\lnsim r\}$.

Vice versa, let in the extended real line $\sup\{r\in\mathbb{R}|r\lnsim X\}=x=\inf\{r\in\mathbb{R}|X\lnsim r\}$.

If $x=+\infty$, then $\sup\{r\in\mathbb{R}|r\lnsim X\}=+\infty$,
which means that if $y$ is a real number, then there is a real number
$r$ such that $y<r$ and $r\lnsim X$. Per \ref{preorder-to-order-of-reals}
also $y\lnsim X$, which proves that $E(X)=+\infty=x$.

If $x=-\infty$, then $\inf\{r\in\mathbb{R}|X\lnsim r\}=-\infty$,
which means that if $y$ is a real number, then there is a real number
$r$ such that $r<y$ and $X\lnsim r$. Per \ref{preorder-to-order-of-reals}
also $X\lnsim y$, which proves that $E(X)=-\infty=x$.

If $x\in\mathbb{R}$ and $\epsilon$ is a positive real number, then
since $\sup\{r\in\mathbb{R}|r\lnsim X\}=x$, there is a real number
$r$ such that $x-\epsilon<r\lnsim X$, and due to \ref{preorder-to-order-of-reals}
also $x-\epsilon\lnsim X$. Since $x=\inf\{r\in\mathbb{R}|X\lnsim r\}$,
there is a real number $r$ such that $X\lnsim r<x+\epsilon$ , and
due to \ref{preorder-to-order-of-reals} also $X\lnsim x+\epsilon$.
Per \ref{Expectation}, $E(X)=x$.

\subsection{Proof of \ref{Chain-infinite-e}}

Let $E\left(X|C.D\right)\in\left\{ +\infty,-\infty\right\} $ and
$p$ be a positive real number such that $p\lesssim_{D}C$. Per \ref{Conditional-preorder-df},
$pD\lesssim C.D$. Define

\[
s=\begin{cases}
1 & \mathrm{if}\:E(X|C.D)=+\infty\\
-1 & \mathrm{if}\:E(X|C.D)=-\infty
\end{cases}
\]

and

\[
Y=sX
\]

Per \ref{Homogeneity}, $E\left(Y|C.D\right)=sE\left(X|C.D\right)=+\infty$.
Let $y$ be a real number. Define $z=\frac{\max\left(0,y\right)}{p}$.
Then $z\geq0$. Since $E\left(Y|C.D\right)=+\infty$ and $z$ is a
real number, $zC.D\lnsim Y.C.D$. Therefore, $yD\lesssim\max\left(0,y\right)D=zpD\lesssim zC.D\lnsim Y.C.D$.
This proves that $E\left(Y.C|D\right)=+\infty$. Since $X=sY$ and
per \ref{Homogeneity}, $E\left(X.C|D\right)=sE\left(Y.C|D\right)=E\left(X|C.D\right)$.

\subsection{Proof of \ref{Chain}}

Let $x=E\left(X|C.D\right)$ be a real number, $c=P\left(C|D\right)$
be a real number and let $\epsilon$ be a positive real number. Define
$\delta=\frac{\epsilon}{1+\left|x\right|}$. Then $\delta$ is a positive
real number and since $x=E\left(X|C.D\right)$, $-\delta C.D\lnsim X.C.D-xC.D\lnsim\delta C.D$.
In the natural order of events $C.D\le D$ and per \ref{preorder-to-natural-order},
$C.D\lesssim D$. Therefore, $-\delta D\lnsim X.C.D-xC.D\lnsim\delta D$.
Since $c=P\left(C|D\right)$, $-\delta D\lnsim C.D-cD\lnsim\delta D$,
implying that $-\left|x\right|\delta D\lesssim xC.D-xcD\lesssim\left|x\right|\delta D$.
Summing inequalities, we get that $-\left(1+\left|x\right|\right)\delta D\lnsim X.C.D-xcD\lnsim\left(1+\left|x\right|\right)\delta D$.
Since $\left(1+\left|x\right|\right)\delta D=\epsilon D$, also $-\epsilon D\lnsim X.C.D-xcD\lnsim\epsilon D$,
which proves that $E\left(X.C|D\right)=xc$.

The only case when the formula $E\left(X|C.D\right).P\left(C|D\right)$
makes sense and at least one of $E\left(X|C.D\right)$, $P\left(C|D\right)$
is not a real number, is the case when $E\left(X|C.D\right)\in\{-\infty,+\infty\}$
and $c=P\left(C|D\right)$ is a positive real number. In this case,
note that $0<\frac{c}{2}<c$, use \ref{Preorder-consistency} to obtain
that $\frac{c}{2}\lnsim_{D}C$ and use \ref{Chain-infinite-e} proven
above to get that $E\left(X.C|D\right)=E\left(X|C.D\right)=E\left(X|C.D\right).P\left(C|D\right)$.

\subsection{Proof of \ref{Chain-zero-e}}

Let $\epsilon$ be a positive real number. If $E\left(X|C.D\right)=0$,
then per \ref{cond-expect-df}, $-\epsilon C.D\lnsim X.C.D\lnsim\epsilon C.D$.
In the natural order of events $C.D\le D$ and per \ref{preorder-to-natural-order},
$C.D\lesssim D$. Therefore, $-\epsilon D\lesssim-\epsilon C.D\lnsim X.C.D\lnsim\epsilon C.D\lesssim\epsilon D$,
which proves that $E\left(X.C|D\right)=0$.

\subsection{Proof of \ref{Chain-zero-p}}

Let $P\left(C|D\right)=0$ and $p$ be a real number such that $-pC.D\lesssim X.C.D\lesssim pC.D$.
Define $q=\max\left(1,p\right)$. Then both $q>0$ and $q\geq p$.
Therefore, $pC.D\lesssim qC.D$. Let $\epsilon$ be a positive real
number. Define $\delta=\frac{\epsilon}{q}$. Since $P\left(C|D\right)=0$,
$-\delta D\lnsim C.D\lnsim\delta D$, i.e. also $qC.D\lnsim q\delta D$.
Therefore, $-\epsilon D=-q\delta D\lnsim-qC.D\lesssim-pC.D\lesssim X.C.D\lesssim pC.D\lesssim qC.D\lnsim q\delta D=\epsilon D$,
which proves that $E\left(X.C|D\right)=0$.

\subsection{Proof of \ref{Cond-form-inf-e}}

Let $E\left(X.C|D\right)\in\{-\infty,+\infty\}$. Define
\[
s=\begin{cases}
1 & \mathrm{if\:}E(X.C|D)=+\text{\ensuremath{\infty}}\\
-1 & \mathrm{if\:}E(X.C|D)=-\text{\ensuremath{\infty}}
\end{cases}
\]
and 
\[
Y=sX
\]
Per \ref{Homogeneity}, $E\left(Y.C|D\right)=s.E\left(X.C|D\right)=+\infty$.
Let $y$ be a real number. Define $z=\max\left(0,y\right)$. Then
both $z\geq0$ and $z\geq y$. Since $E\left(Y.C|D\right)=+\infty$,
$zD\lnsim Y.C.D$. In the natural order of events $C.D\le D$ and
per \ref{preorder-to-natural-order}, $C.D\lesssim D$. Therefore,
$yC.D\lesssim zC.D\lesssim zD\lnsim Y.C.D$, which proves that $E\left(Y|C.D\right)=+\infty$.
Since $X=sY$ and per \ref{Homogeneity}, $E\left(X|C.D\right)=sE\left(Y|C.D\right)=E\left(X.C|D\right)$.

\subsection{Proof of \ref{Cond-form}}

Let $x=E\left(X.C|D\right)$ be a real number, $c=P\left(C|D\right)$
be a nonzero real number. Per \ref{Minimal-and-maximal-pr}, $c>0$.
Let $\epsilon$ be a positive real number. Define $\delta=\frac{c}{2.\left(1+\frac{\left|x\right|}{c}\right)}.\min(1,\epsilon)$.
Then $0<\delta\le\frac{c}{2}$. Equality $c=P\left(C|D\right)$, the
fact that $\delta$ is positive and \ref{Expectation} give $-\delta D\lnsim C.D\lyxmathsym{\textendash}cD\lnsim\delta D$.
Therefore, $\frac{c}{2}D=cD-\frac{c}{2}D\lesssim cD-\delta D\lnsim C.D$,
i.e. $D\lnsim\frac{2}{c}C.D$. By \ref{plaus-preor-multiplicative-property},
$-\delta\frac{\left|x\right|}{c}D\lesssim-\frac{x}{c}C.D+xD\lesssim\delta\frac{\left|x\right|}{c}D$.
Equality $x=E\left(X.C|D\right)$ and \ref{Expectation} give $-\delta D\lnsim X.C.D-xD\lnsim\delta D$.
Summing the inequalities, we get that $-\delta\left(1+\frac{\left|x\right|}{c}\right)D\lnsim X.C.D-\frac{x}{c}C.D\lnsim\delta\left(1+\frac{\left|x\right|}{c}\right)D$.
Also, $\delta\left(1+\frac{\left|x\right|}{c}\right)D\lnsim\delta\left(1+\frac{\left|x\right|}{c}\right)\frac{2}{c}C.D=\min\left(1,\epsilon\right)C.D\lesssim\epsilon C.D$.
Combining these inequalities, we get that $-\epsilon C.D\lnsim-\delta\left(1+\frac{\left|x\right|}{c}\right)D\lnsim X.C.D-\frac{x}{c}C.D\lnsim\delta\left(1+\frac{\left|x\right|}{c}\right)D\lnsim\epsilon D$,
which proves that $E\left(X|C.D\right)=\frac{x}{c}$.

The only remaining case when the expression $\frac{E\left(X.C|D\right)}{P(C|D)}$
makes sense is the case when $E\left(X.C|D\right)\in\left\{ -\infty,+\infty\right\} $
and $P(C|D)$ is a positive real number. In this case, use the equality
$E\left(X|C.D\right)=E\left(X.C|D\right)$ from \ref{Cond-form-inf-e}
proven above and the fact that $E\left(X.C|D\right)=\frac{E\left(X.C|D\right)}{P(C|D)}$
to finally obtain $E\left(X|C.D\right)=\frac{E\left(X.C|D\right)}{P(C|D)}$.

\subsection{Proof of \ref{Cond-form-zero-e}}

Let $E\left(X.C|D\right)=0$ and let $p$ be a positive real number
such that $pD\lesssim C.D$. Let $\epsilon$ be a positive real number.
Define $\delta=\epsilon p$. Since $E\left(X.C|D\right)=0$ and $\delta>0$,
per \ref{Expectation} $-\delta D\lnsim X.C.D\lnsim\delta D$, implying
that $-\epsilon C.D\lesssim-\epsilon pD=-\delta D\lnsim X.C.D\lnsim\delta D=\epsilon pD\lesssim\epsilon C.D$,
which proves that $E\left(X|C.D\right)=0$.

\subsection{Proof of \ref{Cond-form-zero-p}}

Let $P\left(C|D\right)=0$ and let $p$ be a real number such that
$1\lesssim_{D}pX.C$. Defining $Y=pX$ we get $1\lesssim_{D}Y.C$
and per \ref{cond-expect-df}, $D\lesssim Y.C.D$. Let $y$ be a real
number. Define $\epsilon=\frac{1}{\max\left(1,y\right)}$. Then $\epsilon$
is a positive real number and since $P\left(C|D\right)=0$, $-\epsilon D\lnsim C.D\lnsim\epsilon D$.
Therefore, $\frac{1}{\epsilon}C.D\lnsim D$ and $yC.D\lesssim\max\left(1,y\right)C.D=\frac{1}{\epsilon}C.D\lnsim D\lesssim Y.C.D$.
This proves that $E\left(Y|C.D\right)=+\infty$. It also proves that
$0\lnsim Y$, implying that $Y\ne0$ and since $Y=pX$ also $p\neq0$
and $X=\frac{1}{p}Y$. Per \ref{Homogeneity}, $E\left(X|C.D\right)=\frac{1}{p}E\left(Y|C.D\right)=\frac{+\infty}{p}$.

\subsection{Proof of \ref{P-form-inf-e(x|cd)}}

Let $E\left(X|C.D\right)\in\left\{ -\infty,+\infty\right\} $ and
let $p$ be a real number such that $-p\lesssim_{D}X.C\lesssim_{D}p$.
Per \ref{Conditional-preorder-df}, $-pD\lesssim X.C.D\lesssim pD$.
Define 
\[
q=\max\left(1,p\right)
\]
\[
s=\begin{cases}
1 & \mathrm{if}\:E(X|C.D)=+\text{\ensuremath{\infty}}\\
-1 & \mathrm{if}\:E(X|C.D)=-\text{\ensuremath{\infty}}
\end{cases}\enskip\enskip
\]
and 
\[
Y=sX
\]
Then $q>0$ and $q\geq p$. Therefore, $-qD\lesssim-pD\lesssim X.C.D\lesssim pD\lesssim qD$.
Since $Y=sX$ and $s\in\left\{ -1,1\right\} $, also $-qD\lesssim Y.C.D\lesssim qD$.
Per \ref{Homogeneity}, $E\left(Y|C.D\right)=s.E\left(X|C.D\right)=+\infty$.
Let $\epsilon$ be a positive real number. If we define $y=\frac{q}{\epsilon}$,
then $y>0$. Since $E\left(Y|C.D\right)=+\infty$, $yC.D\lnsim Y.C.D$.
Combining inequalities, we get that $yC.D\lnsim qD$. Therefore, $C.D\lnsim\frac{q}{y}D=\epsilon D$.
Per \ref{preorder-to-natural-order}, $0\lesssim C.D$. Combining
inequalities, we get that $-\epsilon D\lnsim-C.D\lesssim0\lesssim C.D\lnsim\epsilon D$,
proving that $P\left(C|D\right)=0$.

\subsection{Proof of \ref{P-form}}

Let $v=E\left(X.C|D\right)$ be a real number and let $x=E\left(X|C.D\right)$
be a nonzero real number. Let $\epsilon$ be a positive real number.
Define $\delta=\frac{\left|x\right|}{2}\epsilon$. Then $\delta$
is a positive real number and since $v=E\left(X.C|D\right)$, $-\delta D\lnsim X.C.D-vD\lnsim\delta D$.
Since $x=E\left(X|C.D\right)$, $-\delta C.D\lnsim X.C.D-xC.D\lnsim\delta C.D$.
In the natural order of events, $C.D\le D$ and per \ref{preorder-to-natural-order}
$C.D\lesssim D$. Therefore, $-\delta D\lesssim-\delta C.D\lnsim-X.C.D+xC.D\lnsim\delta C.D\lesssim\delta D$.
Summing inequalities, we get $-2\delta D\lnsim xC.D-vD\lnsim2\delta D$,
i.e. $-\epsilon D=-\frac{2}{\left|x\right|}\delta D\lnsim C.D-\frac{v}{x}D\lnsim\frac{2}{\left|x\right|}\delta D=\epsilon D$.
This proves that $E\left(C|D\right)=\frac{v}{x}$.

If $E\left(X.C|D\right)\in\left\{ -\infty,+\infty\right\} $, then
per \ref{Cond-form-inf-e}, $E\left(X|C.D\right)=E\left(X.C|D\right)$,
i.e. the expression $\frac{E\left(X.C|D\right)}{E\left(X|C.D\right)}$
does not make sense.

Let $v=E\left(X.C|D\right)$ be a real number and let $E\left(X|C.D\right)\in\left\{ -\infty,+\infty\right\} $.
Define $p=1+\left|v\right|$. Per \ref{Conditional-expectation},
$-1D\lnsim X.C.D-vD\lnsim1D$, i.e. $-pD=-\left(1+\left|v\right|\right)D\lesssim-\left(1-v\right)D\lnsim X.C.D\lnsim\left(1+v\right)D\lesssim\left(1+\left|v\right|\right)D=pD$.
This demonstrates that the assumptions of \ref{P-form-inf-e(x|cd)}
hold. Therefore, $P\left(C|D\right)=0=\frac{v}{E\left(X|C.D\right)}$.

\subsection{Proof of \ref{P-form-zero-e(xc|d)}}

Let $E\left(X.C|D\right)=0$ and let $p$ be a real number such that
$1\lesssim_{C.D}pX$. Define $Y=pX$. Per \ref{Homogeneity}, $E\left(Y.C|D\right)=pE\left(X.C|D\right)=0$.
Also, $1\lesssim_{C.D}Y$. Per \ref{Conditional-preorder-df}, $C.D\lesssim Y.C.D$.
Let $\epsilon$ be a positive real number. Since $E\left(Y.C|D\right)=0$,
$-\epsilon D\lnsim Y.C.D\lnsim\epsilon D$. Combining inequalities,
$C.D\lnsim\epsilon D$. Per \ref{plausible-preord-plausible-property},
$0\lesssim C.D$. Combining inequalities, we get $-\epsilon D\lnsim-C.D\lesssim0\lesssim C.D\lnsim\epsilon D$,
proving that $P\left(C|D\right)=0$.

\subsection{Proof of \ref{Kolmogorovian-coherence}}

Let $PV$ be a function from $\mathcal{F}$ to $\mathbb{R}_{\geq0}$
that satisfies \ref{Kolmogorovian-pv}. Let $n\geq0$, $m\geq1$ be
integers, $q_{1},\ldots,q_{n}$ be nonnegative real numbers, $r_{1},\ldots,r_{m}$,
$s_{1},\ldots,s_{m}$ be real numbers, $C_{1},\ldots,C_{n}$ be events,
$A_{1},\ldots,A_{m}$ be elements of $\mathcal{F}$ and let $r_{j}\left(PV\left(A_{j}\right)+s_{j}\right)>0$
for every $j\in\left\{ 1,\ldots,m\right\} $. 

Let $\mathcal{A}$ be a Boolean algebra generated by $A_{1},\ldots,A_{m}$
and $\mathcal{B}$ be a Boolean algebra generated by $A_{1},\ldots,A_{m},C_{1},\ldots,C_{n}$.
Then $\mathcal{A}\subseteq\mathcal{B}\subseteq\mathcal{F}$. Since
both $\mathcal{A}$ and $\mathcal{B}$ are finitely generated, they
are finite and atomic \cite{key-9}. Let $\mathrm{at}\left(\mathcal{A}\right)$
be the set of atoms of $\mathcal{A}$ and $\mathrm{at}\left(\mathcal{B}\right)$
be the set of atoms of $\mathcal{B}$. Let $\mathrm{span}\left(\mathcal{B}\right)$
be the linear span of $\mathcal{B}$. Let $G\in\mathrm{at}\left(\mathcal{B}\right)$.
If $X\in\mathrm{span}\left(\mathcal{B}\right)$, then there is a unique
real number $r$ such that $X.G=rG$. This allows us to define a function
$\varphi$ from $\mathrm{span}\left(\mathcal{B}\right)\times\mathrm{at}\left(\mathcal{B}\right)$
to $\mathbb{R}$ such that $X.G=\varphi\left(X,G\right)G$. The properties
of the function $\varphi$ are: 
\begin{itemize}
\item If $B\in\mathcal{B}$, then $\varphi\left(B,G\right)\in\{0,1\}$,
implying that $\varphi\left(B,G\right).\varphi\left(B,G\right)=\varphi\left(B,G\right)$.
(idempotence)
\item If $X,Y\in\mathrm{span}\left(\mathcal{B}\right)$ and $X.Y=0$, then
$\varphi\left(X,G\right).\varphi\left(Y,G\right)=\varphi\left(X.Y,G\right)=0$.
(orthogonality)
\item If $X,Y\in\mathrm{span}\left(\mathcal{B}\right)$, then $\varphi\left(X+Y,G\right)=\varphi\left(X,G\right)+\varphi\left(Y,G\right)$.
(additivity)
\item If $r$ is a real number and $X\in\mathrm{span}\left(\mathcal{B}\right)$,
then $\varphi\left(rX,G\right)=r\varphi\left(X,G\right)$. (homogeneity)
\end{itemize}
Define a function $\nu$ from $\mathcal{C}$ to $\mathbb{R}$ so that
for $B\in\mathcal{B}$, $\nu\left(B\right)=\sum_{G\in\mathrm{at}\left(\mathcal{B}\right)}\varphi\left(B,G\right)$.
The properties of the function $\nu$ are
\begin{itemize}
\item If $0\neq B\in\mathcal{B}$, then $\nu\left(B\right)>0$.
\end{itemize}
Define a function $F$ from $\mathrm{span}\left(\mathcal{B}\right)$
to $\mathbb{R}$ so that for $X\in\mathrm{span}\left(\mathcal{B}\right)$,
\[
F\left(X\right)=\sum_{G\in\mathrm{at}\left(\mathcal{B}\right)}\sum_{H\in\mathrm{at}\left(\mathcal{A}\right)}\varphi\left(X,G\right).\varphi\left(H,G\right).\frac{PV\left(H\right)}{\nu\left(H\right)}
\]
The properties of the function $F$ are
\begin{itemize}
\item If $B\in\mathcal{B}$, then $F\left(B\right)\geq0$. (nonnegativity) 
\item $F$ is additive.
\item $F$ is homogeneous. 
\item $F$ coincides with $PV$ on $\mathcal{A}$.
\end{itemize}
Therefore, $\sum_{j=1}^{m}r_{j}\left(F\left(A_{j}\right)+s_{j}F\left(1\right)\right)=\sum_{j=1}^{m}r_{j}\left(PV\left(A_{j}\right)+s_{j}PV\left(1\right)\right)>0$.
Also, $F\left(\sum_{i=1}^{n}q_{i}C_{i}\right)\geq0$ and $F\left(\sum_{i=1}^{n}q_{i}C_{i}+\sum_{j=1}^{m}r_{j}\left(A_{j}+s_{j}\right)\right)>0$.
Since $F\left(0\right)=0$, this proves that $\sum_{i=1}^{n}q_{i}C_{i}+\sum_{j=1}^{m}r_{j}\left(A_{j}+s_{j}\right)\neq0$.

\subsection{Proof of \ref{Coxian-sum-rule}}

Let $PV$ be a function from $\mathcal{F}\times\mathcal{F}_{0}$ to
$\mathbb{R}_{\geq0}$ that satisfies definition \ref{Coxian-pv}.
Let $A,B$ be elements of $\mathcal{F}$ such that $A.B=0$ and $C\in\mathcal{F}_{0}$.
Then $A+B=1-(1-A).(1-B)$ is an element of $\mathcal{F}$.

If $(1-B).C=0$, then $B.C=C$, i.e. $PV(B|C)=PV(B.C|C)=PV(C|C)=1$.
Also $A.C=A.B.C=0$ implying $PV(A|C)=PV(A.C|C)=0$. Finally, $PV(A+B|C)=PV((A+B).C|C)=PV(A.C+B.C|C)=PV(C|C)=1$.
Therefore, $PV(A+B|C)=PV(A|C)+PV(B|C)$.

If $(1-B).C\neq0$, then $PV(A+B|C)=1-PV((1-A).(1-B)|C)=1-PV((1-A)|(1-B).C).PV(1-B|C)=1-[1-PV(A|(1-B).C)].PV(1-B|C)=1-PV(1-B|C)+PV(A|(1-B).C).PV(1-B|C)=PV(B|C)+PV(A.(1-B)|C)=PV(B|C)+PV(A-A.B|C)=PV(B|C)+PV(A|C)$.

\subsection{Proof of \ref{Coxian-coherence}}

Let $PV$ be a function from $\mathcal{F}\times\mathcal{F}_{0}$ to
$\mathbb{R}_{\geq0}$ that satisfies definition \ref{Coxian-pv}.
Let $n\geq0,m\geq1$ be integers, $q_{1},\ldots,q_{n}$ be nonnegative
real numbers, $r_{1},\ldots,r_{m}$, $s_{1},\ldots,s_{m}$ be real
numbers, $C_{1},\ldots,C_{n}$ be events, $A_{1},\ldots,A_{m}$ be
elements of $\mathcal{F}$, $D_{1},\ldots,D_{m}$ be elements of $\mathcal{F}_{0}$
and for every $j\in{1,\ldots,m}$, $r_{j}(PV(A_{j}|D_{j})+s_{j})>0$.

Let $D=\bigvee_{j=1}^{m}D_{j}$. Since $PV$ is subadditive, $PV(\bigvee_{j=1}^{m}D_{j}|D)\leq\sum_{j=1}^{m}PV(D_{j}|D)$.
Since $0<PV(D|D)$, there is a $k\in\{1,\ldots,m\}$ such that $0<PV(D_{k}|D)$.

Let $j\in{1,\ldots,m}$. Due to Bayes' rule and nonnegativity, 
\[
r_{j}(PV(A_{j}.D_{j}|D)+s_{j}PV(D_{j}|D))=r_{j}(PV(A_{j}|D_{j})+s_{j})PV(D_{j}|D)\geq0
\]
If $k\in\{1,\ldots,m\}$ is such that $0<PV(D_{k}|D)$, then $r_{k}(PV(A_{k}D_{k}|D)+s_{k}PV(D_{k}|D))>0$.
Since such $k\in\{1,\ldots,m\}$ exists, $\sum_{j=1}^{m}r_{j}(PV(A_{j}D_{j}|D)+s_{j}PV(D_{j}|D))>0$.

Let $\mathcal{A}$ be the Boolean algebra generated by $A_{1},\ldots,A_{m},D_{1},\ldots,D_{m}$.
Since $\mathcal{A}$ is finitely generated, it is finite and atomic
\cite{key-9}. Let $\mathrm{at}\left(\mathcal{A}\right)$ be the set
of atoms of $\mathcal{A}$.

Let $\mathcal{B}$ be the Boolean algebra generated by $A_{1},\ldots,A_{m},D_{1},\ldots,D_{m},C_{1},\ldots,C_{n}$.
Since $\mathcal{B}$ is finitely generated, it is finite and atomic
\cite{key-9}. Let $\mathrm{at}\left(\mathcal{B}\right)$ be the set
of atoms of $\mathcal{B}$ and $\mathrm{span}\left(\mathcal{B}\right)$
be the linear span of $\mathcal{B}$.

Define a function $\varphi$ from $\mathrm{span}\left(\mathcal{B}\right)\times\mathrm{at}\left(\mathcal{B}\right)$
to $\mathbb{R}$ so that if $X\in\mathrm{span}\left(\mathcal{B}\right)$
and $G\in\mathrm{at}\left(\mathcal{B}\right)$, then $X.G=\varphi(X,G).G$.
Also define a function $\nu$ from $\mathcal{B}$ to $\mathbb{R}$
so that if $H\in\mathcal{B}$, then $\nu(H)=\sum_{G\in\mathrm{at}\left(\mathcal{B}\right)}\varphi(H,G)$.
Finally, define a function $F$ from $\mathrm{span}\left(\mathcal{B}\right)\times\{D\}$
to $\mathbb{R}$ so that if $X\in\mathrm{span}\left(\mathcal{B}\right)$,
then 
\[
F(X|D)=\sum_{G\in\mathrm{at}\left(\mathcal{B}\right)}\sum_{H\in\mathrm{at}\left(\mathcal{A}\right)}\varphi(X,G).\varphi(H,G).\frac{PV(H|D)}{\nu(H)}
\]
The reader can verify that $F$ is homogeneous and additive on $\mathrm{span}\left(\mathcal{B}\right)\times\{D\}$,
nonnegative on $\mathcal{B}\times\{D\}$ and that it coincides with
$PV$ on $\mathcal{A}\times\{D\}$. Therefore, 
\begin{align*}
F(\sum_{j=1}^{m}r_{j}(A_{j}+s_{j}).D_{j}|D) & =\sum_{j=1}^{m}r_{j}(F(A_{j}.D_{j}|D)+s_{j}F(D_{j}|D))\\
 & =\sum_{j=1}^{m}r_{j}(PV(A_{j}.D_{j}|D)+s_{j}PV(D_{j}|D))>0
\end{align*}
Also, $F(\sum_{i=1}^{n}q_{i}C_{i}|D)\geq0$. Therefore, $F\left(\sum_{i=1}^{n}q_{i}C_{i}+\sum_{j=1}^{m}r_{j}(A_{j}+s_{j}).D_{j}|D\right)>0$.
Since $F(0|D)=0$, this proves that $\sum_{i=1}^{n}q_{i}C_{i}+\sum_{j=1}^{m}r_{j}(A_{j}+s_{j}).D_{j}\neq0$.

\subsection{Proof of \ref{Dupre-Tiplerian-coherence}}

Let $PV$ be a partial function from $\mathcal{T}\times\mathcal{C}$
to $\mathbb{R}$ that satisfies definition \ref{Dupre-Tiplerian-pv}.
Let $n\geq0,m\geq1$ be integers, $q_{1},\ldots,q_{n}$ be nonnegative
real numbers, $r_{1},\ldots,r_{m}$, $s_{1},\ldots,s_{m}$ be real
numbers, $C_{1},\ldots,C_{n}$, $D_{1},\ldots,D_{m}$ be events, $X_{1},\ldots,X_{m}$
be random quantities and for every $j\in{1,\ldots,m}$, $r_{j}(PV(X_{j}|D_{j})+s_{j})>0$.

Due to additivity of $PV$, if $C\in\mathcal{C}$, then $PV(0|C)=PV(0|C)+PV(0|C)$,
which implies that $PV(0|C)=0$.

Let $D=\bigvee_{j=1}^{m}D_{j}$. Due to subadditivity, $PV(D|D)\leq\sum_{j=1}^{m}PV(D_{j}|D)$.
Since $0<PV(D|D)$, there is a $k\in\{1,\ldots,m\}$ such that $0<PV(D_{k}|D)$.

Let $j\in{1,\ldots,m}$. Due to Bayes' rule and nonnegativity, 
\[
r_{j}(PV(X_{j}.D_{j}|D)+s_{j}PV(D_{j}|D))=r_{j}(PV(X_{j}|D_{j})+s_{j}).PV(D_{j}|D)\geq0
\]
 If $k\in\{1,\ldots,m\}$ is such that $0<PV(D_{k}|D)$, then $r_{k}(PV(X_{k}.D_{k}|D)+s_{k}PV(D_{k}|D))>0$.
Since such $k\in\{1,\ldots,m\}$ exists, $\sum_{j=1}^{m}r_{j}(PV(X_{j}.D_{j}|D)+s_{j}PV(D_{j}|D))>0$.
Additivity and homogeneity of $PV$ imply that $PV\left(\sum_{j=1}^{m}r_{j}(X_{j}+s_{j}).D_{j}|D\right)>0$.
Moreover, $PV(\sum_{i=1}^{n}q_{i}C_{i}|D)\geq0$, i.e. also $PV\left(\sum_{i=1}^{n}q_{i}C_{i}+\sum_{j=1}^{m}r_{j}(X_{j}+s_{j}).D_{j}|D\right)>0$.
Since $PV(0|D)=0$, this proves that $\sum_{i=1}^{n}q_{i}C_{i}+\sum_{j=1}^{m}r_{j}(X_{j}+s_{j}).D_{j}\neq0$.

\subsection{Proof of \ref{Coherent functions}}

We start by proving that $1\Rightarrow2$.

Let $PV$ be a coherent partial function from $\mathcal{T}\times\mathcal{E}\left(\mathcal{T}\right)$
to $\overline{\mathbb{R}}$.

Define the relation $\lnsim$ so that $0\lnsim X$ if $X=\sum_{i=1}^{n}p_{i}C_{i}+\sum_{j=1}^{m}r_{j}\left(X_{j}+s_{j}\right).D_{j}$
for some nonnegative integers $n$, $m$, positive real numbers $p_{1},\ldots,p_{n}$,
real numbers $r_{1},\ldots,r_{m}$, $s_{1},\ldots,s_{m}$, nonzero
events $C_{1},\ldots,C_{n}$, events $D_{1},\ldots,D_{m}$ and random
quantities $X_{1},\ldots,X_{m}$, such that at least one of $n$,
$m$ is nonzero and $r_{j}\left(PV\left(X_{j}|D_{j}\right)+s_{j}\right)>0$
for every $j\in\left\{ 1,\ldots,m\right\} $. Define $X\lnsim Y$
if $0\lnsim Y-X$. According to this definition, for every nonzero
event $C$ holds that $0\lnsim C$. This guarantees both regularity
and plausible property \ref{strict-p-o-pl-property}. The antireflexivity
of $\lnsim$ is a consequence of the coherence of $PV$ and \ref{pos-comb-nz-events}.
We leave the task to verify that $\lnsim$ has the remaining fundamental
properties of a plausible strict partial order listed in \ref{strict-p-o-fund}
as an exercise to the reader.

Take $\lesssim$ as a plausible preorder having $\lnsim$ as its strict
part. By \ref{strict-p-o-sufficiency}, such a plausible preorder
exists.

If $PV\left(X|C\right)=+\infty$ and $y$ is a real number, then $PV\left(X|C\right)-y=+\infty>0$,
i.e. $yC\lnsim X.C$, proving that $E\left(X|C\right)=+\infty$.

If $PV\left(X|C\right)=-\infty$ and $y$ is a real number, then $\left(-1\right).\left(PV\left(X|C\right)-y\right)=+\infty>0$,
i.e. $X.C\lnsim yC$, proving that $E\left(X|C\right)=-\infty$.

If $r=PV\left(X|C\right)$ is a real number and $\epsilon$ is a positive
real number, then $PV\left(X|C\right)-r+\epsilon>0$, i.e., $-\epsilon C\lnsim X.C-rC$,
and $\left(-1\right).\left(PV\left(X|C\right)-r-\epsilon\right)>0$,
i.e. $X.C-rC\lnsim\epsilon C$, proving that $E\left(X|C\right)=r$.
This completes the proof that the expectation naturally induced by
$\lesssim$ extends the function $PV$.

The implication $2\Rightarrow3$ is trivial.

We finish by proving that $3\Rightarrow1$.

Let $\lesssim$ be a plausible preorder and $E$ be the conditional
expectation induced by $\lesssim$. Let $n\geq0$, $m\geq1$ be integers,
$q_{1},\ldots,q_{n}$ be nonnegative real numbers, $r_{1},\ldots,r_{m}$,
$s_{1},\ldots,s_{m}$ be real numbers, $C_{1},\ldots,C_{n}$, $D_{1},\ldots,D_{m}$
be events, $X_{1},\ldots,X_{m}$ be random quantities and for every
$j\in\{1,\ldots,m\}$, $r_{j}\left(E\left(X_{j}|D_{j}\right)+s_{j}\right)>0$.
By \ref{pl-preord-df}, $0\lesssim\sum_{i=1}^{n}q_{i}C_{i}$. By real
additivity and homogeneity, $E\left(r_{j}\left(X_{j}+s_{j}\right)|D_{j}\right)=r_{j}\left(E\left(X_{j}|D_{j}\right)+s_{j}\right)>0$.
By preorder consistency, if $E\left(r_{j}\left(X_{j}+s_{j}\right)|D_{j}\right)>0$,
then $0\lnsim r_{j}\left(X_{j}+s_{j}\right).D_{j}$, yielding $0\lnsim\sum_{i=1}^{n}q_{i}C_{i}+\sum_{j=1}^{m}r_{j}\left(X_{j}+s_{j}\right).D_{j}$.
Due to \ref{strict-p-o-Antireflexivity}, $0\neq\sum_{i=1}^{n}q_{i}C_{i}+\sum_{j=1}^{m}r_{j}\left(X_{j}+s_{j}\right).D_{j}$.
\end{document}